\documentstyle{article}
\begin{document}

\noindent
Ukrainian Mathematical Journal, Vol. 46, No. 12, 1994, 1907-1914.

\bigskip

\centerline {\large \bf MINIMAL HANDLE DECOMPOSITION OF}

\medskip
\centerline {\large \bf SMOOTH SIMPLY CONNECTED 5-MANIFOLDS}

\bigskip
\centerline {\large \bf A. O. Prishlyak}

\bigskip

\centerline {Kiev University, Ukraine}

\medskip

e-mail: prish@mechmat.univ.kiev.ua


\bigskip
A theorem on the existence of the unique minimal topologic handle decomposition
of differentiable simply connected five-dimensional manifolds is proved. For a
decomposition of this sort, the number of handles of each index is given.

\bigskip
Assume that $W^n$ is a smooth compact simply connected manifold with boundary
$\partial W= V_0 U\ V_1$ where  $V_0$ and $V_1$ are connected components of the
boundary. We say that the sequence of inclusions $W_0 \subset W_1\subset  ..$.
$\subset  W$ is a handle decomposition of $W$ beginning from $V_0\times  I$ if

\medskip
(i) $W_0$ is a collar $V_0 \times  I;$

\medskip
(ii) $W_{i+1}$ results from $W_i$ after pasting a handle of index $k$ to
$\partial W_i\backslash V_0$, $i.e.$, $W_{i+1}=
(W_i.U\ D^k\times D^{n-k})\backslash ~$, where ~ denotes the identification of
the points  $\partial D^k\times D^{n-k} =S^{k-1} \times D^{n-k}$ with the images
of these points in $\partial W_i\backslash V_0$ for some inclusion  $f: S^{k-1}
\times D^{n-k} \rightarrow  \partial W_i\backslash V_0;$

(iii) for any $x\in  W$, there exists a neighborhood that intersects only finitely
miny sets $W_i\backslash W_i-_1$ [1],

\medskip
Let $N_\lambda $ be a number of handles of index $\lambda $ in a handle
decomposition of W. This decomposition is called minimal if there is no other
handle decomposition of $W$ with the entire number of handles less than
$\Sigma \ N_\lambda $ and with the number of handles of index $\lambda $ at
most $N_\lambda $ for each $\lambda $. We say that a minimal handle
decomposition is unique if any other minimal handle decomposition has the same
number of handles of index $\lambda $ for each $\lambda $

The existence of the unique smooth minimal handle decomposition for smooth
simply connected manifolds was proved by Smale in [2] (in the case of simply
connected boundaries for $n \leq  6)$, by Barden in [3] (in the case of a
closed manifold for $n = 5)$, and by Shkol'nikov in [4] (in the case
of manifolds with standard boundaries for $n = 5$).

We say that the triple  $( C, G, d )$  is a crossed module if $C$ and $G$
are groups such that $G$ acts on $C$ from the left and $d: C \rightarrow  G$ is
a homomorphism such that $c_1 + c_2 - c_1 = d(c_1)c_2$,  $d(gc) = g(d(c))
g-^1$. Sharko has shown that if $n \leq  6$, then, on $W^n$, there exists a
unique minimal handle decomposition without handles of indices 0, 1, $n - 1$,
$n$, in which the number of handles of index 2 is equal to the minimal number
of generators of the group $\pi _2 (W^n\backslash V_1)$  considered as a
$\pi _1(V_1)$-module [5].

Freedman has proved that if $W^5$ is a simply connected five-dimensional
$h$-cobordism, then $W^5$ is homeomorphic to $V_0 \times  I$. By using the
method of proving this assertion, one can show that there exists a unique
topological handle decomposition of a smooth simply connected manifold $W^5$
with simply connected boundary components $V_0$ and $V_1$. If the boundary
components are not simply connected, then the transversal spheres are absent
and, hence, this method cannot be applied to the belt spheres of two-handles,
which are necessary for constructing the Kasson handles.

\medskip

{\bf Lemma.}Assume that $W$ is a simply connected compact five-dimensional
manifold with the boundary $\partial W= V_0 \bigcup V_1$, where $V_i$,
$i = 0$, 1, are
connected components of the boundary. Then there exists a minimal system of
generators $\alpha _1, \alpha _2, ..., \alpha _k$ of the group $\pi _2(W$,
$V_1)$ considered as a crossed $\pi _1(V_1)$ - module such that, under the
action of the Gurevich homomorphism $\pi _2(W$, $V_1) \rightarrow  H_2(W$,
$V_1)$, the generators $\alpha _1, \alpha _2, ..., \alpha _s$, $s \leq
k$, should turn into a minimal system of generators of the group $H_2(W$,
$V_1)$ and the generators $\alpha _{s+1}, ... , \alpha _k$ should lie
in the kernel of this homomorphism.

{\bf Proof.} Consider a commutative diagram consisting of the segments of exact
homotopic and exact homo-logical sequences of the pair ( $W$, $V_1)$

$$
\matrix {\ \pi_{2}(W) & \to & \pi_{2}(W,V_{1}) & \to & \pi_{1}(V_{1}) & \to
& 0 \cr
\downarrow & & \downarrow & &\downarrow & & \cr
H_{2}(W) & \to & H_{2}(W,V_{1}) & \to & H_{1}(V_{1}) & \to & 0 \cr }
$$
\medskip
Since

$$\pi _2(W) = H_2(W),\  \pi _1(V_1) /[\pi _1(V_1), \pi _1(V_1)] = H_1(V_1)$$
\medskip
\noindent
$\pi _2(W,V_1)\rightarrow  H_2(W,V_1)$ is an epimorphism. Let
$\beta _1,..., \beta _k$ be an arbitrary minimal system of
generators of the group $\pi _2(W,V_1)$ considered as a crossed
$\pi _1(V_1)$-module. Denote by $\omega _i$ . the image of $\beta _i$.
under the homomorphism $\pi _2(W,V_1)\rightarrow  H_2(W,V_1)$. Then $\omega _i$
form a system of generators of the group $H_2(W,V_1)$. The following operations
are called elementary transformations:

\medskip

(i) the change of $\omega _i$ by $\omega _i \pm \omega _j.$, where $i \neq j;$

\medskip
(ii) the multiplication of a generator $\omega _i$ by -1;

\medskip
(iii) the renumbering of generators.

\medskip
\noindent
It is clear that, under elementary transformations, a system of generators
turns into a system of generators. Let us show that, for the given system, the
elementary transformations over the generators $\omega _i$ enable one to select
a system with the minimal number of generators $\gamma _1$,
$\gamma _2,...,\gamma _s.$

Assume that the system $\gamma _1$, $\gamma _2,...,\gamma _S$ with the minimal
number of generators corresponds to the decomposition of a group into a direct
sum of cyclic groups $(\gamma _i$ are generators of the direct summands). Then

\medskip
\noindent
$$\gamma _{i} = \sum_{j=1}^{k} a_{ij}\omega _{j}, \
\omega _{j}=\sum_{i=1}^{m} b_{ij} \gamma _{i}.$$

\medskip
\noindent
If $\gamma _i$. has order $p_i$ then the coefficient $b_{ij}$ is taken modulo
$p_i$. For a fixed $i$, we consider the coefficients $b_{i1}$, $b_{i2},...$,
$b_{im}$. If $\gamma _i$ is a free generator, then $b_{i1}$, $b_{i2},...$,
$b_{im}$ have the greatest common divisor equal to one. Therefore, after
elementary transformations, $\omega _j$ can be reduced to the form 0, ... , 0,
1, 0,... , 0 . By substituting the new $\omega _j$ in the relation

$$\gamma_{i} = \sum_{j=1}^{k} a_{ij}\omega_{j}$$

\medskip
\noindent
we verify that the coefficient $a_{ij}$ is equal to one. In other words,

$$\gamma _i = \omega _i + \sum_{i\neq j} a_{ij}\omega _j$$

\medskip
\noindent
and, hence, after elementary transformations, $\omega _i$ can be replaced by
$\gamma _i$. If $\gamma _i$ has order $p_i$, then, by elementary transformation
of $\omega _j.$, we can reduce the coefficients $b_{i1}$, $b_{i2},...$,
$b_{im}$ to the form 0,...,0, $q_i$, 0,..., 0, where $q_i$ and $p_i$ are
relatively prime numbers and, hence, $q_i\gamma _i$ is also a generator. After
multiplying the relation

$$\gamma _i = \sum_{j=1}^{k}{a_{ij}\omega _j},$$

\medskip
\noindent
by $q_i$ we obtain that, on the basis of the elementary transformations, one
can replace $\omega _i$ by $q_i\gamma _i$. Thus, a system with the minimal
number of generators can be selected from the system of generators $\omega _i.$

Let us perform similar operations over the generators $\beta _1,...,
\beta _k$. Namely, if, when reducing to a minimal system of
generators, $\omega _i$. in $H_2(W,V_1)$ is replaced by $\omega _i.\pm
\omega _j$, then the corresponding generator $\beta _i$ in $\pi _2(W$,
$V_1)$ is replaced by $\beta _i\pm  \beta _j$. Thus, a new
system of generators $\alpha _1,\alpha _2$, ... , $\alpha _k$ will remain a
minimal one for the group $\pi _2(W$, $V_1)$, the generators
$\alpha _1,\alpha _2$, ... , $\alpha _S$ will turn into the minimal system of
generators of the group $H_2(W$, $V_1)$, and the generators $\alpha _{m+1}$,
... , $\alpha _k$ will lie in the commutator $[\pi _2(W$, $V_1)$, $\pi _2(W$,
$V_1)].$

\medskip
{\bf Theorem.} Assume that $W$ is a smooth simply connected compact
five-dimensional manifold with the boundary $\partial W= V_0 U\ V_1$. Denote by
$V_0$ the boundary component for which

$$\mu (\pi _2(W, V_0)) -\mu (H_2(W, V_0)) \leq  \mu (\pi _2(W,V_1))
-\mu (H_2(W, V_1)).$$

\medskip
\noindent
Then, on $W$, there exists the unique minimal topological handle decomposition
without handles of indices 0, 1,4, and 5, but with $\pi _2(W$, $V_1)$ handles
of index 2 and $\mu (\pi _2(W$, $V_1))-\mu (\pi _2(W$, $V_1)+\mu (\pi _2(W$,
$V_1)$ handles of index 3. Here, $\mu (H)$ is a minimal number of generators of
the group $H$ and $\pi _2(W$, $V_i)$, $i=l,2$, is regarded as a crossed
$\pi _1(V_i)${\it -}module.

\medskip
{\bf Proof.} Let $W$ be a smooth manifold satisfying the conditions of the
theorem. Denote

$$k := \mu (\pi _2(W, V_0)),$$

$$s :=\mu (H_2(W, V_0)),$$

$$m := \mu (\pi _2(W, V_1)),$$

$$n :=\mu (H_2(W, V_1)).$$

\medskip
Consider a handle decomposition of the manifold, beginning from $V_0\times I$.
Since $\pi (W,V_i)=0$, $i=0,1$, one can find a handle decomposition without
handles of indices 0, 1,4, and 5 [6, 7]. Let $h^\alpha _1,..$. , $h^\alpha _k$
be tube neighborhoods of the disks that form a minimal system of generators of
the group  $\pi _2(W,V_0)$  regarded as a $\pi _1(V_0$, $x)$-module, and let
$H^\alpha _1,...$, $H^\alpha _m$ be tube neighborhoods of the disks that form a
minimal system of generators of the group $\pi _2(W$, $V_1)$ regarded as a
$\pi _1(V_0$, $y)$-module.  In addition, let these neighborhoods be such that
the core of $h^\alpha _1,...$, $h^\alpha _s$ and  $H^\alpha _1,...$,
$H^\alpha _n$ are minimal systems of generators of the groups  $H_2(W$, $V_0)$
and  $H_2(W,V_1)$, respectively. Below, the described neighborhoods play the roles of 2-
and 3-handles. Consider a handle decomposition of $W$ which contains the
handles $h^\alpha _1,...$, $h^\alpha _k$ and $H^\alpha _1,...$, $H^\alpha _m$,
and the other 2- and 3- handles lie in$cl(W\setminus (\cup h_{i}\cup
H_{i}^{\alpha}))$ . To this handle decomposition, we add
the following $k + m$ pairs of mutually canceling 2-and 3-handles: the
2-handles $(h_{1}^{\alpha},\cdots ,h_{k^{\alpha}}, h_{1}^{\gamma},\cdots ,
h_{m}^{\gamma}$ and the corresponding 3-handles $H_{1}^{\beta},\cdots ,H_{k}^{\beta}, H_{1}^{\gamma},\cdots
,H_{m}^{\gamma}$. After sliding the handle $h^{\beta }_i$ over $h^{\alpha }_i$, we
replace the handles $h^{\beta }_i$ by $h_i= h^{\alpha }_i
+h^{\beta }_i$ and the handles $H^{\gamma }_i$ by $H_i =
H^{\alpha }_i+ H^{\gamma }_i$. Then the index is defined as

$$\lambda
(h_{i},H_{j}^{\gamma}) =\delta_{i}^{j} , \lambda
(h_{i}^{\gamma},H_{j})=\delta_{i}^{j}, (\delta_{i}^{i}=1; \
\delta_{i}^{j}=0,i\ne j)$$

\medskip
and the matrix A consisting of the indices of intersection of 2- and 3-handles
has the form

$$
\matrix{\ & \matrix{\ \quad H_{1}^{\beta}\ & \ldots & H_{k}^{\beta} &\ H_{1}
& \ldots & H_{m} & & \ast & & } \cr
\cr
\matrix{\ h_{1} \cr
\vdots \cr
h_{k} \cr
\cr
h_{1}^{\gamma} \cr
\vdots \cr
h_{k}^{\gamma} \cr
\cr
\ast \cr
} & \left | \matrix{\ \quad 1 \quad & & 0 & & & & & & \cr
& \ddots & & & 0 & & & \ast & \cr
0 & &\ 1\ & & & & & & \cr
\cr
& & &\quad 1 & & & & & \cr
& 0 & & & \ddots & & & 0 & \cr
& & & & & \quad 1 & & & \cr
\cr
& 0 & & &\ast & & &\ \ast \ & \cr
} \right | } $$

\noindent
By adding the handles, we reduce this matrix to the form

$$
\left | \matrix{\ E & 0 & 0 \cr
0 & E & 0 \cr
0 & 0 & B \cr
} \right |
$$

\noindent
in this case, the core of $h_i$ and $H_i$ as $h^\alpha _i$ and $H^\alpha _i$,
should be minimal systems of generators of the groups $\pi _2(W,V_o)$ and
$\pi _2(W,V_1)$ regarded as crossed modules.

Below, we assume that this handle decomposition has no collar $V_0\times I$
between 2- and 3-handles.

We set $W^{1}=cl(W\setminus (\cup h_{i}\cup H_{j}\cup h_{p}^{\gamma} \cup
H_{q}^{\beta})), \ \partial W=V^1_0 U\ V^1_1$, where $V^1_i$ are connected components of
the boundary. Then $W^1$ consists of the remaining 2- and 3-handles. For this
manifold, the homology groups $H_2(W^1,V^1_0)$ and $H_2\{W^1,V^1_1)$ are
specified by the matrix $B$ and coincide with the corresponding homology groups
of the manifold $W$ because the handles $h_i$ and $H^{\beta }_i$,
$h^\gamma _j$ and $H_j$ have an index of intersection equal to one and do not
contribute to homology groups.

On the basis of the van Kampen theorem, one can show that $\pi _1(W^1)
=\pi _1(V^1_0 ) = \pi _1(V^1_1) = 0$. By adding handles on the manifold $W^1$,
we reduce the matrix $B$ to the diagonal form such that, on the diagonal, the
numbers $\pm  1$ should correspond only to the handles   $h^{\alpha }_1,..$. ,
$h^{\alpha }_s$ and  $H^{\alpha }_1,...$, $H^{\alpha }_n$ whose cores are minimal
systems of generators of the groups  $H_2(W, V_0)$  and   $H_2(W, V_1)$.
We set  $W^{2}= cl(W^{1}\setminus (\cup h_{i}^{\alpha} \cup
H_{j}^{\alpha})), 1\le i\le s, 1\le j\le n$. Then the matrix
consisting of the indices of intersection of 2- and 3-handles on $W^2$ is the
identity matrix and we have the homologies  $H_i(W^2, V^2_0)=0, 1 \le i \le 4$,
 where  $V^2_0$ is  a  component  of  the  boundary $\partial W^2$.   By
the Whitehead theorem, $W^2$ is an $h$-cobordism, and it follows from the
Freedman theorem that $W^2$ is homeomorphic to $V^2_0 \times I [8].$

By construction, the cores of the handles $h_1,... , h_s$ and $H_1,...
,H_n$ specify the same elements of the
groups $H_2(W, V_0)$ and $H_2(W, V_1)$ as $h^{\alpha }_1,... , h^{\alpha }_s$
and $H^{\alpha }_1,... , H^{\alpha }_n.$ Let us change their places- In
this case, if we cut out all the handles from $W$ except $h_1,... , h_s$,
$H^{\beta }_1,... , H^{\beta }_s$, and $h^{\gamma }_1,... ,
h^{\gamma }_s,$ $H_1,..., H_s$, we obtain a manifold with simply connected
boundaries, which is a simply connected $h$-cobordism and, by the Freedman
theorem [8], has the structure of topological product.

Thus, the matrix A can be reduced to the form

$$
\matrix{\ & \matrix{\ \quad H_{1}^{\omega } \ \ldots \ \ H_{n}^{\omega } & H_{n+1}
\ \ldots \ \ H_{m} &  H_{s+1}^{\gamma } \ \ldots \ \ H_{k}^{\gamma } & } \cr
\cr
\matrix{\ h_{1}^{\omega } \cr
\vdots \cr
h_{s}^{\omega } \cr
\cr
h_{s+1} \cr
\vdots \cr
h_{k} \cr
\cr
h_{n+1}^{\gamma} \cr
\vdots \cr
h_{m}^{\gamma}
 \cr
} & \left | \matrix{\ \quad n_1  & & 0 & & & & & & \cr
& \ddots & & & 0 & & & 0 & \cr
0 & & n_t\ & & & & & & \cr
\cr
& & & & & &\quad 1 & & \cr
& 0 & & & 0 & &  & \ddots   \cr
& & & & & & & & \ 1 \quad  \cr
\cr
& & &\quad 1 & & & & & \cr
& 0 & & & \ddots & & & 0 & \cr
& & & & & \ 1 & & & \cr
} \right | } $$

\bigskip

The convention on notation of the boundary components $V_i, \ i=0,l$, implies
that $k-s \leq m-n$. By construction, the cores of the handles
$h^{\alpha }_1,... , h^{\alpha }_s,$ $ h_{s+1}, ... , h_k$ are generators of
$\pi _2(W, V_0)$,  $H^{\alpha }_1,... , H^{\alpha }_n$, $H_{n+1}, ... , H_m$
are generators of $\pi _2(W, V_1)$, and $h_{n+1}^{\gamma}\cdots ,h_{m}^{\gamma}$ and
$H_{s+1}^{\gamma},\cdots ,H_{k}^{\gamma}$  are trivial elements
of these groups. Replace $ H_{s+i}^{\beta}$ by $H_{n+i}^{\omega} = H_{s+i}^{\beta}+H_{n+i},$
$ 1\le i\le m-n, $ and $h_{n+i}^{\gamma} $ by $h_{s+i}^{\omega}=h_{n+i}^{\gamma}-
h_{s+i},$ $ 1\le i\le m-n$. Then the cores of the handles specify the same elements
of the group $\pi _2(W$, $V_0)$ as $h_{s+i}$ . We set
$h^{\omega }_i=h^{\alpha }_i$, $1\le i\le s$, and $H^{\omega }_i=H^{\alpha }_i$,
$1\le i\le n$. Then the cores $h^{\omega }_i$ and $H^{\omega }_i$ define
minimal systems of generators of the groups $\pi _2(W, V_0)$ and $\pi _2(W,V_1)$,
and the matrix A can be reduced to the form

\medskip
$$
\matrix{\ & \matrix{\quad \quad H_{1}^{\omega }\quad \quad \quad \ldots \ \quad \quad \quad
H_{m}^{\omega } \ H_{s+m-n+1}^{\beta  } \ldots   H_{k}^{\beta } & } \cr
\cr
\matrix{\ h_{1}^{\omega } \cr
\ \cr
\ \cr
\cr
\ \cr
\vdots \cr
\ \cr
\ \cr
\ \cr
\ \cr
h_{m}^{\gamma}
 \cr
} & \left | \matrix{\ \quad n_1  & & 0 & & & & & & \cr
& \ddots & & & 0 & & & 0 & \cr
0 & & n_t\ & & & & & & \cr
\cr
& &    &\quad 1 & & & & & \cr
& 0 & &      & \ddots  & & & 0 & \cr
& & & &  & \ 1  & & & \cr
\cr
& & & & & & \quad 1 &   & \cr
& 0 & & &  0 & & & \ddots     & \cr
& & & & & & & & \ 1  \quad  \cr
} \right | } $$
\bigskip

Thus, we have constructed a handle decomposition with $k$ handles of index 2
and $k - s + n$ handles of index 3.

Let us prove that this handle decomposition is minimal. Consider an arbitrary
handle decomposition of $W$, beginning from $V_0\times I$. Let $W_2$ be a
submanifold containing all the handles of indices at most 2 and only these
handles. Let us show that the number $q$ of critical handles of index 2 in this
handle decomposition is at least $k$. The mapping $\pi _2(W_2, V_1)\rightarrow
\pi _2(W, V_0)$ is an epimorphism because any element of
the group $\pi _2(W, V_0)$  is realized by an inclusion of a two-dimensional
disk and, by general reasons, can be isotoped into the manifold $W_2$. By
applying the Whitehead theorem of the structure of the second relative
homotopic group, one can show that the cores of handles of index 2 define the
generators of the group $\pi _2(W, V_0)$ and, hence, their number is at least
$k$. Let handles $h_1, ... , h_k$ and $H_1, ... , H_k$ from the given
decomposition specify the minimal systems of generators of the groups
$\pi _2(W, V_0)$  and $\pi _2(W, V_1)$, respectively. We set $X=\bigcup h_i$, $Y=
\bigcup H_i$. Consider the manifold $Z=W\backslash (X \bigcup V)$ with the boundary
$\partial Z= V^1_0 \bigcup V^1_1$ where $V^1_0=\partial X\backslash V_0$,
$V^1_1=\partial Y\backslash V_1$. By using the van Kampen-Seifert theorem, one
can show that

$$\pi _1(Z)= \pi _1(Y)= \pi _1(X)= \pi _1(V^1_0)= \pi _1(V^1_1)= 0.$$

\medskip
\noindent
The epimorphism $\pi _2(X$, $V_0) \rightarrow  \pi _2(W$, $V_0)$ induces the
epimorphism $\pi _2(X) \rightarrow  \pi _2(W)$ in the exact homotopic sequence
of the pair $(W, X)$. Hence, the groups $\pi _2(W$, $X)$ and $H_2(W,X)$ are
trivial. In this case, by the Poincare theorem on excision and duality,

$$h_2(w,x) = h_2(Y \cup Z, V^1_0) = h_3(Y \cup Z, V_1) = 0.$$
\medskip
\noindent
Since $W$ is obtained from $Y \bigcup Z$ by pasting 3-handles, we have

\medskip

$$H_3(X, V^1_0) = H_3(W, Y \cup  Z) = kZ,    H_2(X, V^1_0) = 0.$$

\medskip

Similarly,

\medskip

$$\pi _2(W,Y) = H_2(W,Y) =H_2(X \cup Z,Vi^1) = H_3(X \cup Z,V_0) = 0,$$

\medskip

$$H_3(Y,V^1_1) = H_3(W,X \cup Z) = mZ,  H_2(Y, V^1_1)) = 0.$$

\bigskip
\noindent
Consider the segment of the exact homological sequence of the triple
$(X \bigcup Z,$ $X,$ $V_0)$

$$
\matrix {\ 0 \to H_{3}(X\cup Z,X) \to H_{3}(X\cup Z,X) \to H_{3}(X\cup Z,X)
\to H_{2}(X\cup Z,X) \to 0 \cr
\Vert \qquad \ \ \qquad \qquad \Vert \qquad \qquad \ \ \qquad \Vert \ \ \
\qquad \qquad \qquad \Vert \cr
H_{3}(Z,V_{0}^{1}) \ \qquad \qquad kZ \qquad \qquad \qquad \lambda Z \qquad
\ \qquad H_{2}(Z,V_{0}^{1}) \cr
}
$$

\medskip
\noindent
where $\lambda $ can be found from the exact homological sequence of the triple
 $(W, X \bigcup Z,$ $V_0)$

\medskip
$$0\rightarrow  H_3(W, V_0) \rightarrow  H_3(W, X \cup Z) \rightarrow
H_2(X \cup Z,V_0) \rightarrow  H_2(W, V_0) \rightarrow  0,$$

\medskip
$$\lambda = \mu (H_2(X \cup Z,V_0)) = m + s - n.$$

\medskip
Then

\medskip
$$\mu (H_3(Z,V_0)) - \mu (H_2(Z,V^1_0)) =k-\lambda = k-m+n-s.$$

\medskip
\noindent
and $\mu (H_3(Z,V_0)) \leq  k-\lambda = k-m+n-s$. It follows from the
Morse inequality that, on the manifold $Z$, there exist at least $k-m+n-s$
handles of index 3. Thus, on the manifold $W$, there exist at least $k+n-s$
handles of index 3, whence the handle decomposition constructed in order to
prove the theorem is minimal.

{\bf Corollary 1}. Assume that $W$ is a smooth simply connected compact
five-dimensional manifold with the connected boundary $\partial \ W =V_0$ .
Then, on $W$, there exists a unique minimal topological handle decomposition
without handles of indices 0, 1, 4, and 5, but with $\mu (\pi _2(W$, $V_0))$
handles of index 2 and $\mu (\pi _2(W$, $V_0)) - \mu (H_2(W$, $V_0 )) +
\mu (H_2(W$, $V_1))$ handles of index 3. Here, $\mu (H)$ is the minimal number
of generators of the group $H$ and $\pi _2(W,V_i)$, $i=l,2$, is regarded as a
crossed $\pi _1(V_i)$ -module.

{\bf Corollary 2}. Assume that $W$ is a compact contractible five-dimensional
manifold with boundary. Then, on $W$, there exists a minimal topological handle
decomposition with one handle of index 5, without handles of indices 0, 1,
and4, but with $\mu (\pi _2(W$, $\partial W))$ handles of indices 2 and 3.
Here$,\mu (H)$ is the minimal number of generators of the group $H$ and
$\pi _2(W,\partial W)$ is regarded as a crossed $\pi _1(\partial W)$-module.

\medskip
This research was partially supported by the International Scientific
Foundation, grant No, V6F000.

\bigskip

\centerline {\bf REFERENCES}

\medskip
1. $R$, Kirby and $L$, Sibenman, ``Foundational essay on topological manifolds,
smoothing and triangulations," Ann, Math Stud,, 88 (1977).

2. S. Smale, ``Generalized Poincare's conjecture in dimensions greater than
four, " Ann. Math.  74, No, 2, 391-406 (1961).

3. D. Barden, ``Simply connected five-manifolds," Ann, Math, 82,365-385 (1965).

4. Yu, A. Shkol'nikov, ``Handle decomposition of simply connected
five-manifolds. III," Ukr,Mat.Zh,, 46, No.7,935-940 (1994).

5. V. $V$, Sharko, Functions on Manifolds [in Russian], Naukova Dumka, Kiev
(1990).

6. C. P. Rourke and B. J. Sanderson, Introduction to Piecewise-Linear Topology,
Springer, Berlin (1972).

7. J. Milnor, Lectures on the $h$-Cobordism Theorem, Princeton University
Press, Princeton (1965).

8. M. Freedman, ``The topology of four-dimensional manifolds, " J. Different.
Geom., 17, 357-453 (1982).
\end{document}